\documentclass[12pt]{amsart}
\usepackage{amssymb,latexsym,amsmath,amscd}
\usepackage[latin1]{inputenc}
\newcommand{\PP}{\mathbb{P}}
\newcommand{\ZZ}{\mathbb{Z}}

\newcommand{\OO}{\mathcal{O}}

\newcommand{\LL}{\mathcal{L}}
\newcommand{\MM}{\mathcal{M}}
\newcommand{\NN}{\mathcal{N}}
\newcommand{\EE}{\mathcal{E}}
\newcommand{\FF}{\mathcal{F}}
\newcommand{\GG}{\mathcal{G}}

\newcommand{\JJ}{\mathcal{J}}

\newcommand{\lra}{\longrightarrow}
\newcommand{\ra}{\rightarrow}

\newcommand{\ffi}{\varphi}

\newcommand{\coker}{\operatorname{Coker}\,}
\newcommand{\srel}{\stackrel}
\renewcommand{\deg}{\operatorname{deg}}

\newcommand{\rg}{\operatorname{rk}}

\newcommand{\depth}{\operatorname{depth}}

\theoremstyle{definition}
\newtheorem{defi}{Definition}[section]
\theoremstyle{plain}
\newtheorem{lema}[defi]{Lemma}
\newtheorem{teo}[defi]{Theorem}

\newtheorem{prop}[defi]{Proposition}
\newtheorem{cor}[defi]{Corollary}
\newtheorem{rk}[defi]{Remark}

\newtheorem{ex}[defi]{Example}

\newtheorem{teo-def}[defi]{Theorem/Definition}
\newenvironment{proofteo}{\noindent {\textit{Proof of Theorem \ref{stable}, continued.}}}{\quad \hfill $\Box$}

\title{ACM bundles on cubic surfaces}
\author{Marta Casanellas}
\address{Departament de Matemàtica Aplicada I. ETSEIB. Universitat Polit\`ecnica de Catalunya. Avinguda Diagonal 647. 08028 Barcelona. Spain.}
\email{marta.casanellas@upc.edu}
\author{Robin Hartshorne}
\address{Department of Mathematics. Evans Hall. University of California. Berkeley, CA, 94720-3840. USA}
\email{robin@math.berkeley.edu}
\thanks{Research of the first author partially supported by Ministerio de Educaci\'on y Ciencia MTM2006-E14234-C02-02.}
\begin{document}
\maketitle
\begin{abstract}
%In this paper we prove that, for every $r \geq 2$, the moduli
%space $M^s_X(r;c_1,c_2)$ of rank $r$ stable vector bundles with
%Chern classes $c_1=rH$ and $c_2=\frac{1}{2}(3r^2-r)$ on a
%nonsingular cubic surface  $X \subset \mathbb{P}^3$ contains a
%nonempty smooth open subset formed by ACM bundles. The bundles we
%consider for this study are extremal for the number of generators
%of the corresponding module (these are known as Ulrich bundles),
%so  we also prove the existence of indecomposable Ulrich bundles
%of arbitrarily high rank on $X$.
In this paper we prove that, for every $r \geq 2$, the moduli
space $M^s_X(r;c_1,c_2)$ of rank $r$ stable vector bundles with
Chern classes $c_1=rH$ and $c_2=\frac{1}{2}(3r^2-r)$ on a
nonsingular cubic surface $X \subset \mathbb{P}^3$ contains a
nonempty smooth open subset formed by ACM bundles, i.e. vector
bundles with no intermediate cohomology. The bundles we consider
for this study are extremal for the number of generators of the
corresponding module (these are known as Ulrich bundles), so  we
also prove the existence of indecomposable Ulrich bundles of
arbitrarily high rank on $X$.
%In this paper we exhibit families of of rank $r$ stable ACM bundles (i.e. without intermediate cohomology) on a nonsingular cubic surface $X \subset \mathbb{P}^3$ for each $r \geq 2$. As a consequence we prove that these bundles form a smooth open subset of the moduli space of rank $r$ stable vector bundles with Chern classes $c_1=rH$ and $c_2=\frac{1}{2}(3r^2-r)$, $M^s_X(r;c_1,c_2)$. The bundles we consider for this study are extremal for the number of generators of the corresponding module (these are known as Ulrich bundles).
\end{abstract}

\section{Introduction}

One has known since the work of Mumford that to have a reasonable parameter space for families of algebraic vector bundles on an algebraic variety $X$, it is necessary to impose a condition of stability. The existence of a coarse moduli space of stable vector bundles has been known since the work of Maruyama \cite{Maruyama77}, Gieseker \cite{Gieseker} and Simpson \cite{Simpson}. The moduli of vector bundles of all ranks on algebraic curves has been studied extensively. Bundles of rank 2 on higher dimensional varieties have often played an important role --think of the work \cite{Mukai} of Mukai on the moduli of stable rank 2 bundles on a K3 surface. By means of the Serre correspondence, rank 2 bundles on any variety are related to the study of codimension 2 subvarieties (see for example Beauville \cite{Beauvillecubic}). However, little is known about bundles of higher rank.

A particular class of vector bundles that have been studied in recent years is the class of arithmetically Cohen-Macaulay (ACM) bundles. A vector bundle $\EE$ on a projective variety $X$ is ACM if all its intermediate cohomology groups $H^i(X, \EE(m))$ are zero for $0 <i<\dim X$ and all $m \in \ZZ$. These bundles correspond to Maximal Cohen-Macaulay modules (MCM) over the associated graded ring. In the algebraic context, MCM modules have been extensively studied (see for example the book of Yoshino \cite{Yoshino}), as they reflect relevant properties of the corresponding ring. There has also been recent work on ACM bundles of small rank on particular varieties such as Fano 3-folds, quartic threefolds and Grassmann varieties (see \cite{AM} and the references therein).

From the point of view of representation theory, one can often
divide varieties (or rings) into three classes depending on the
behavior of the category of ACM bundles (or MCM modules). A
variety is of \textit{finite representation type} if there are
only finitely many indecomposable ACM bundles (up to twist). The
projective varieties of finite representation type have been
completely classified into a small list (see \cite{Yoshino}). For
example, on a nonsingular quadric hypersurface $X$ in $\PP^n$,
there are only one (or two) indecomposable ACM bundles (not
counting the structure sheaf $\OO_X$), depending on the parity of
$n$. A variety is of \textit{tame representation type} if, for each rank $r$, the indecomposable ACM
bundles of rank $r$ form a finite number of families of dimension at most one. For example,
according to the work of Atiyah \cite{Atiyah}, on an elliptic
curve, for each rank and degree there is a single family of
indecomposable bundles, parameterized by the curve itself. A
variety is of \textit{wild representation type} if
there exist $n$-dimensional families of non-isomorphic
indecomposable ACM bundles for arbitrarily large $n$. Since all vector bundles on projective
curves are ACM, we see that a nonsingular curve in projective
space is of finite, tame, or wild representation type according as
its genus is 0, 1, or $\geq 2$. Drozd and Greuel
\cite{DrozdGreuel92} have shown that  the category of MCM modules
over the complete local ring of a reduced curve singularity is
either finite, tame or wild. However, we cannot expect such a
trichotomy in general, because for example a quadric cone in
$\PP^3$ has an infinite discrete set of indecomposable ACM bundles
of rank 2 \cite[6.1]{CH}.
% If the variety is neither finite nor tame, then it is called of \textit{wild representation type}, and one expects families of arbitrarily high dimension of ACM bundles of arbitrarily high rank. For example, on curves of genus $\geq 2$ there are families of indecomposable stable bundles of arbitrarily large dimension. See the work \cite{Drozd} of Y. Drozd on the trichotomy finite-tame-wild.

In this paper we exhibit families of stable ACM bundles on a nonsingular cubic surface in $\PP^3$, of arbitrarily high rank and dimension. %This is the first case we know of specific families of stable ACM bundles of higher rank on varieties of dimension grater than one of wild representation type.
Thus the cubic is of wild representation type. As far as we know, these are the first examples of indecomposable ACM bundles of arbitrarily high rank on any varieties except curves.

Our bundles have another particular characteristic: they are extremal for the number of generators of the corresponding module. In the case of MCM modules, this phenomenon was discovered by Ulrich \cite{Ulrich}. He showed that for MCM modules over a local ring, there is a bound on the number of generators of the  module in terms of the multiplicity and the rank. Since then, modules with the maximum number of generators have been called Ulrich modules, and correspondingly Ulrich bundles in the geometric case. Even the existence of Ulrich bundles on  projective varieties is not known in general, though in our case of a cubic surface in $\PP^3$ it is easy, because the ideal sheaf of a twisted cubic curve is one. Also rank two ACM bundles on the cubic surface have been classified by Faenzi \cite{Faenzi}, and in his classification there are stable Ulrich bundles with $c_1=2H$, $c_2=5$ that form a family of dimension 5.
%Here we prove that there exist Ulrich bundles of any rank on any nonsingular cubic surface $X$ and that they form families of dimension $r^2+1$ showing that $X$ is of wild type.

Here is our main result, \eqref{mainteo}.

\textbf{Theorem.}\textit{ Let $X$ be a nonsingular cubic surface
in $\PP^3$ over an algebraically closed field $k$. Then for every
$r\geq 2$ there are stable Ulrich bundles of rank $r$ with
$c_1=rH$ and $c_2=\frac{1}{2}(3r^2-r)$.
%The corresponding moduli space is smooth of dimension $r^2+1$.
They form a smooth open subset of an irreducible component of dimension $r^2+1$ of the moduli space $M^s_X(r;c_1,c_2)$ of rank $r$ stable vector bundles with Chern classes $c_1$, $c_2$ on $X$.}

% Although it looks like Ulrich bundles are very special (and they are indeed), we prove that the bundles of the previous theorem form an open set of the moduli space of stable vector bundles of rank $r$:
%
% \textbf{Corollary.}\textit{
% The moduli space $M^s_X(r;c_1,c_2)$ of rank $r$ stable vector bundles on $X$ with Chern classes $c_1 =rH$, $c_2=\frac{1}{2}(3r^2-r)$ has an irreducible component of dimension $r^2+1$ containing a smooth open subset that corresponds to stable Ulrich bundles.}

The paper is organized as follows. In section 2 we recall properties of ACM sheaves including the theory of matrix factorizations due to Eisenbud. In section 3 we recall the definition and some properties of Ulrich sheaves including a new proof of the bound on the number of generators of the graded module.

In section 4 we give the construction of our Ulrich bundles over the cubic surface, based on an earlier work of the first author on minimal free resolutions of sets of points in general position, which in turns depends on earlier work of the second author on Gorenstein liaison of zero-schemes on the cubic surface. In section 5 we show that our bundles are stable, and compute the dimension of the coarse moduli space.

In the last section 6 we show that even though our surface is of wild representation type, our bundles are not so wild after all.
We show that these bundles are stably equivalent to \textit{layered} ACM bundles \eqref{cor_Ulrichlayered}, where layered means a successive extension of rank 1 ACM bundles (corresponding to ACM curves in the surface). In particular our bundles are direct summands of successive extensions of ACM line bundles.
% bundles with a filtration whose quotients are either rank 2  layered ACM sheaves or syzygies of such. Layered means an extension of rank 1 ACM sheaves corresponding to ACM curves on the surface. In other words, our sheaves can be related to the ideal sheaves of ACM curves on the surface, by the operations of successive extension, syzygies and stable equivalence.

\section{Generalities on ACM sheaves}
Throughout this section $X$  denotes an integral hypersurface in
$\PP^n$ of dimension $\geq 2$ defined by a homogeneous polynomial
$f$ of degree $d$, $R$ is the polynomial ring $k[x_0, \dots, x_n]$
over an algebraically closed field and $R_X$ is the coordinate
ring of $X$. If $\FF$ is an $\OO_X$-module, $H^i_{\ast}(\FF)$
denotes $\oplus_{l \in \ZZ}H^0(X, \FF(l))$ and $h^i(\FF)$ is the
dimension of $H^i(\FF)$. The minimal number of generators of $\FF$
will refer to the minimal number of generators of
$H^0_{\ast}(\FF)$ and will be denoted by $\mu(\FF)$. The dual of a
module $M$ is $M^{\vee}=Hom_{R_X}(M,R_X)$, and analogously for
sheaves. For a torsion-free coherent sheaf $\FF$, $\rg(\FF)$ will
denote its rank. When $\FF \cong \OO_X(a_1) \oplus \dots \oplus
\OO_X(a_r)$ we will say that $\FF$ is a \textit{dissoci\'e} sheaf.

We recall that a subscheme $X \subset \PP^n$ is
\textit{arithmetically Cohen-Macaulay}, \textit{ACM} for short,
(respectively \textit{arithmetically Gorenstein},\textit{ AG} for
short) if its homogeneous coordinate ring $R_X$ is a
$^{\ast}$local graded Cohen-Macaulay ring (resp. Gorenstein ring),
in the sense of \cite{BroadmannSharp}. For example, when $X
\subset \PP^n$ is a hypersurface of degree $d$, then it is arithmetically Gorenstein.

We will deal both with maximal Cohen-Macaulay $R_X$-modules and
arithmetically Cohen-Macaulay coherent sheaves over $X$. We recall
the definition here.
\begin{defi}
A coherent sheaf ${\mathcal E}$ on $X$ is an {ACM \em sheaf} if it
is locally Cohen-Macaulay on $X$ and $H_*^i({\mathcal E})= 0$ for
$0 < i < \dim X$. A graded $R_X$-module $E$ is a maximal
Cohen-Macaulay module (MCM from now on) if $\depth E=\dim E=\dim R_X$.
\end{defi}

There is a one-to-one correspondence between ACM sheaves on
arithmetically Cohen-Macaulay schemes $X$ and graded MCM
$R_X$-modules sending $\EE$ to $H^0_{\ast}(\EE)$ (see \cite[2.1]{CH}
). When $X$ is non-singular, ACM sheaves are locally free, so
we will be speaking about vector bundles in this case.

\begin{defi} For any coherent sheaf ${\mathcal E}$ on the hypersurface $X$ we
define the {\textit{syzygy sheaf}} ${\mathcal E}^{\sigma}$ to be
the sheafification of the kernel of a minimal free presentation of
$E = H_*^0({\mathcal E})$ over $R_X$.  In other words, let $F_0
\rightarrow E \rightarrow 0$ be a minimal free cover, and sheafify
to get
\[
0 \rightarrow {\mathcal E}^{\sigma} \rightarrow {\mathcal F_0}
\rightarrow {\mathcal E} \rightarrow 0
\]
with ${\mathcal F_0}$ free and $H_*^0({\mathcal F_0}) \rightarrow
H_*^0({\mathcal E})$ surjective.
\end{defi}

When $X$ is a hypersurface defined by a polynomial $f$, Eisenbud
proved in \cite{Eisenbud} that there is a one-to-one
correspondence between MCM $R_X$-modules and matrix
factorizations of $f$. We recall here how it works.

If $E$ is a MCM $R_X$-module, then as an $R$-module it has
projective dimension 1 and therefore it has a minimal free
resolution over $R$:
$$0 \lra F \srel{\varphi}{\lra} G  \lra E \lra 0$$
where $F$ and $G$ are free $R$-modules of the same rank and $\ffi$
is a degree zero morphism (i.e $\ffi$ is a \textit{homogeneous}
matrix). As $E$ is annihilated by $f$, $fG \subseteq im(\varphi)$,
so there exists a morphism $\psi:G \lra F$ so that $\varphi
\psi=fId_{G}$. Then, $\varphi \psi \varphi =f \ffi$ and as $\ffi$
is injective, this implies that $\psi \ffi =f Id_F$.
\begin{defi}[\cite{Eisenbud}]
A \textit{matrix factorization} of $f \in R$ is an ordered pair of
morphisms of free $R$-modules $(\ffi:F \ra G, \psi:G \ra F)$ such
that $\ffi \psi=f\textrm{Id}_G$ and $\psi \ffi=f\textrm{Id}_F$.

If $f$ is a nonzero-divisor, then one equality implies the other.
\end{defi}

We just saw that to each MCM $R_X$-module $E$ we
can associate a matrix factorization $(\ffi,\psi)$ with $E \cong
\coker(\ffi)$. Moreover, if $f$  is irreducible, $\det(\ffi)=f^r$
where $r=\rg(E)$ \cite[5.6]{Eisenbud}. A matrix
factorization $(\ffi:F \ra G, \psi:G \ra F)$ is \textit{reduced}
if $\ffi(F) \subset \mathfrak{m}G$ and $\psi(G) \subset
\mathfrak{m}F$ or, in other words, there is no scalar entry in
$\ffi$ or $\psi$ different from 0 (here $\mathfrak{m}$ denotes the
maximal irrelevant ideal of $R$). The correspondence between
matrix factorizations and MCM modules over hypersurface rings can
be summarized as follows:

\begin{teo}[{\cite[6.3]{Eisenbud}}]\label{eisteo} Let $X \subset \PP^n$ be an integral hypersurface defined by a form $f$.
There are bijections between the sets of
\begin{enumerate}
\item Equivalence classes of reduced matrix factorizations $(\ffi,\psi)$ of $f$ over $R$.
\item Isomorphism classes of non-trivial 2-periodic minimal free resolutions over
$R_X$,
$$\dots \lra F(-d) \srel{\ffi}{\lra} G(-d) \srel{\psi}{\lra} F \srel{\ffi}{\lra} G \lra  M=\coker(\ffi)\lra 0 \,.$$
\item MCM $R_X$-modules $M=\coker(\ffi)$ without free summands.
\end{enumerate}
The correspondence between (i) and (iii) sends $(\ffi,\psi)$ to
$M:=\coker\ffi$. Moreover $\det(\ffi)=f^r$ where
$r=\rg(M)$.
\end{teo}

\begin{lema}\label{propertiesACM} If $X \subset \PP^n$ is an integral hypersurface of degree $d$ and $\EE$ is an ACM sheaf over $X$ corresponding to a matrix factorization $(\ffi,\psi)$ then
\begin{enumerate}
\item[(i)] $\EE^{\sigma}$ is also ACM and has no direct free summands.
\item [(ii)] If $\EE$ is indecomposable, then so is $\EE^{\sigma}$.
\item[(iii)]
If $\EE$ has no direct free summands, then $(\psi, \ffi)$ is a
matrix factorization for $\EE^{\sigma}$.
\item[(iv)] $\EE^{\vee}$ is also ACM and $\mu(\EE)=\mu(\EE^{\vee})$. Furthermore $\EE$ is reflexive.
\item[(v)] If $\EE$ has no direct free summands, then  $(\ffi^{\vee},\psi^{\vee})$ is a matrix factorization for $\EE^{\vee}$.
\item[(vi)] If $\EE$ has no direct free summands, then $\EE^{\sigma\sigma} \cong \EE(-d)$,  $\EE^{\sigma\vee} \cong
\EE^{\vee\sigma}(d)$, $\EE^{\sigma\vee\sigma\vee} \cong \EE$.
\item[(vii)] If $0 \rightarrow {\mathcal E}' \rightarrow {\mathcal E}
\rightarrow {\mathcal E}'' \rightarrow 0$ is an exact sequence of
ACM sheaves, then there is a dissoci\'e sheaf ${\mathcal M}$ and
an exact sequence
\[
0 \rightarrow {\mathcal E}'{}^{\sigma} \rightarrow {\mathcal
E}^{\sigma} \oplus {\mathcal M} \rightarrow {\mathcal
E}''{}^{\sigma} \rightarrow 0.
\]
\end{enumerate}
\end{lema}
\begin{proof}
(i) and (ii) follow from  \cite[4.2]{CDH}. Statement (iii) follows
from \ref{eisteo}. From \cite[2.3]{CH} we have that $\EE^{\vee}$ is
also ACM and reflexive. The assertion saying that  $\mu(\EE)$ is equal to $\mu(\EE^{\vee})$ can be
found in \cite[1.5]{HerzogKuhl} and the the proof goes as follows.
Let $E=H_{\ast}^0(\EE)$, $m=\mu(E)$ and let $0 \lra F \lra G \lra
E \lra 0$ be a minimal free $R$-resolution of $E$, with
$\rg(G)=\rg(F)=m$. Then dualizing we obtain that $m=\mu(Ext^1_R(E,R))$.
On the other hand, applying $Hom(E, \,\,)$ to the exact sequence
$0 \lra R(-d) \srel{f}{\lra} R \lra R_X \lra 0$ we obtain the
exact sequence $$0 \lra Hom_R(E,R_X) \lra Ext^1_R(E,R(-d))
\srel{f}{\lra} Ext^1_R(E,R).$$ As $f$ annihilates $E$, the first
module is isomorphic to $Hom_{R_X}(E,R_X)=E^{\vee}$ and the last
morphism is 0. Therefore $m=\mu(E^{\vee})$ and statement (iv)
follows.

Theorem \ref{eisteo} and (iv) imply item (v).

Statement (vi) is a consequence of (iii) and (v) and of the
2-periodic resolution in Theorem \ref{eisteo}.

The last claim was proved in \cite[4.1]{CDH}.
\end{proof}

\section{Generalities on Ulrich sheaves}

The minimal number of generators $\mu(\EE)$ of an ACM
$\OO_X$-module $\EE$ of positive rank over an integral scheme $X$
of degree $d$ is bounded above.
\begin{teo} Let $X \subset \PP^n$ be an integral subscheme and $\EE$ and ACM sheaf on $X$. Then $\mu(\EE)$ (the minimal number of generators of $H^0_{\ast}(\EE)$)
is bounded by
\begin{equation}
\label{ineq}\mu(\EE) \leq \deg(X)\rg(\EE)
\end{equation}
\end{teo}

\begin{rk}\rm
For Cohen-Macaulay local rings this theorem was proved in
\cite{Ulrich} and \cite{BHU}.
\end{rk}

\begin{proof} Let $m=\dim X$, and choose a finite projection $\pi:
X \lra \PP^m$. Then $\pi_{\ast}(\EE)$ is a coherent
sheaf on $\PP^m$ which is locally Cohen-Macaulay because $\depth
\EE$ can be calculated on $X$ or on $\PP^m$ and is the same. Hence
$\pi_{\ast}(\EE)$ is locally free on $\PP^m$.

Moreover, $\pi_{\ast}(\EE)$ is ACM on $\PP^m$. Indeed, since $\pi$
is a finite morphism, $H^i(X, \EE(l))=H^i(\PP^m,
\pi_{\ast}(\EE)(l))$ for all $i, l$. Then by the Theorem of
Horrocks \cite{Horrocks}, $\pi_{\ast}(\EE)$ is dissoci\'e. Now
$\mu(\EE)$ is the minimal number of generators of
$E=H^0_{\ast}(\EE)$ as an $R_X=k[x_0, \dots, x_n]/I_X$-module.
Inside $R_X$ we have the polynomial ring $S=k[x_0, \dots, x_m]$
coming from $\PP^m$. Since $\pi_{\ast}(\EE)$ is dissoci\'e of rank
$\deg( X) \rg(\EE)$, it is minimally generated by exactly $\deg(X)
\rg(\EE)$ elements. These elements also generate $H^0_{\ast}(\EE)$
over $R_X$ but may no longer be minimal, so $\mu(\EE) \leq \deg
(X) \rg(\EE)$.
\end{proof}

Modules attaining the upper bound were studied by
Ulrich in \cite{Ulrich} and named after him in
\cite{HerzogKuhl} (in \cite{BHU} they were called MGMCM).

\begin{defi} Let $X \subset \PP^n$ be an integral scheme. An ACM sheaf $\EE$ over $X$ (respectively MCM $R_X$-module) is an \textit{Ulrich sheaf}
(resp. \textit{Ulrich module}) if the minimal number of generators
of $\EE$ is $\mu(\EE)=\deg(X)\rg(\EE)$ (resp.
$\mu(E)=\deg(X)rk(E)$).

When $X$ is nonsingular, $\EE$ is locally free, so we call it an
\textit{Ulrich bundle}.
\end{defi}

\begin{defi} Let $\EE$ be a coherent sheaf on $X$ (respectively, a finitely generated graded $R_X$-module). Then we say
that $\EE$ is \textit{normalized} if $H^0(X,\EE)\neq 0$ and $H^0(X,
\EE(-1))=0$ (resp. $E_0 \neq 0$, $E_{-1}=0$).
\end{defi}

\begin{cor}\label{boundh0} Let $X \subset \PP^n$ be an integral subscheme of degree $d$ and $\EE$ and ACM sheaf on
$X$. Assume furthermore that $\EE$ is normalized. Then $h^0 (\EE)
\leq d \rg(\EE)$ and equality implies that $\EE$ is an Ulrich
sheaf.
\end{cor}
\begin{proof} Since $\EE$ is normalized, the elements of
$H^0(\EE)$ are part of a minimal system of generators for $\EE$.
Then $h^0(\EE) \leq \mu(\EE) \leq  d\rg(\EE)$. Equalities imply
that $\EE$ is an Ulrich sheaf.
\end{proof}

Here we state some properties of Ulrich ACM sheaves.
\begin{lema}\label{lemEUlrich} Let $\EE$ be a rank $r$ ACM sheaf over an integral
hypersurface $X \subset \PP^n$ of degree $d \geq 2$. Then the following holds:
\begin{itemize}
\item[a)] If $\EE$ is Ulrich then it has no free summands.
\item[b)] $\EE$ is Ulrich if and only if $\EE^{\vee}$ is Ulrich.

%In this case, if $E$ is normalized and has a linear resolution, $E^{\vee}(-d-1)$ is normalized.
\end{itemize}
 \end{lema}
\begin{proof}
a) If $\EE$ is Ulrich and can be decomposed as $\EE=\FF \oplus
\MM$, where $\MM$ is dissoci\'e, then $\mu(\EE)=\mu(\FF)+\rg(\MM)$
and $\FF$ is  an ACM sheaf of rank $r-\rg(\MM)$. Then inequality
(\ref{ineq}) applied to $\FF$ yields $\mu(\EE) \leq
d(r-\rg(\MM))+\rg(\MM)$. But $\mu(\EE)=dr$ by assumption, so this
is a contradiction if $d > 1$.

b) is a consequence of a), \ref{propertiesACM} (iv) and (v).
\end{proof}

\begin{prop}\label{resolUlrich}
Let $\EE$ be a normalized ACM sheaf of rank $r$ over an integral hypersurface $X \subset \PP^n$ of degree $d \geq 2$.
\begin{itemize}\item[a)] If $\EE$ is Ulrich and is not an extension of lower rank Ulrich bundles then it has a linear minimal free resolution over $\OO_{\PP^n}$:
 $$0 \lra \OO_{\PP^n}(-a-1)^{dr} \srel{\varphi}\lra \OO_{\PP^n}^{dr}(-a) \lra \EE \lra 0.$$
\item[b)] The following are equivalent
\begin{enumerate}
\item $\EE$ has a linear minimal free resolution over $\OO_{\PP^n}$
\item $\EE$ is Ulrich and generated by its global sections
\item $h^0(\EE)=dr$
\end{enumerate}
\end{itemize}
%Moreover, if these conditions hold, then $c_2(\EE) \leq
%\frac{1}{2}(3r^2-r)$ and  we have equality if $\EE$ is orientable.
\end{prop}

%
% \item[(ii)] If $E$ has a minimal linear resolution over $R$, then
% it is an Ulrich module.
% \item[(iii)] Conversely, if $E$ is Ulrich and is not an extension of lower rank Ulrich bundles then it has a linear minimal free resolution over $R$:
% $$0 \lra R(-a-1)^{dr} \srel{\varphi}\lra R^{dr}(-a) \lra E \lra 0.$$
% \item[(iv)] Assume that $\EE$ is normalized. Then $h^0(\EE)=dr$ if and only if it has a minimal linear resolution (so it is Ulrich).
% \item[(v)] Assume that $\EE$ is Ulrich and normalized. Then it is
% generated by global sections if and only if it has a linear
% resolution.

\begin{proof}
We call $E$ the module $H_{*}^0(\EE)$.

a) Let $0 \lra \oplus_{j=1}^{dr} R(b_j)  \srel{\ffi}{\lra}
\oplus_{i=1}^{dr} R(a_i) \lra E \lra 0$ be the minimal free
resolution of $E$ as an $R$-module, with $a_1 \geq \dots \geq a_{dr}$
and $b_1 \geq \dots \geq b_{dr}$, $(\varphi,\psi)$  a matrix
factorization of $f$. After choosing bases, $\ffi$ is a matrix
with entries of degree $u_{i,j}=a_i-b_j$, which decrease from top
to bottom and from right to left. As $\EE$ has no free summands,
$\ffi$ is reduced by Theorem \ref{eisteo}, i.e. any entry of
degree 0 is actually 0. In particular, if $u_{i,i}\leq 0$ then
$\det(\ffi)=0$ which is impossible because $(\ffi,\psi)$ is a
matrix factorization for $f$. Arguing similarly for $E^{\sigma}$
we obtain that $u_{i,i}$ must be strictly smaller than $d$, so $1
\leq u_{i,i} \leq d-1$. Moreover, if $u_{i,i-1} \leq 0$, then
$\ffi$ has the form
$$\ffi=\left(\begin{array}{cc}
A_1 & A \\
0 & A_2
\end{array}\right) \,.
$$
If we write $$\psi=\left(\begin{array}{cc}
B_1 & C \\
D & B_2
\end{array}\right)$$ where $B_i$ has the same size as $A_i$, then from $\ffi\psi=fId$ and $\psi\ffi=fId$ we obtain that $(A_1,B_1)$ and $(A_2,B_2)$ are matrix factorizations of $f$ because $B_1A_1=fId$, $A_2B_2=fId$. Therefore there is an extension sequence $$0 \lra \coker A_1 \lra E \lra \coker A_2 \lra 0.$$ Moreover $\coker A_1$ and $\coker A_2$ are Ulrich because if $A_i$ is an $l_i\times l_i$ matrix, then $\mu(E) \leq \mu(\coker A_1)+\mu(\coker A_2) \leq dl_1+dl_2 =dr$ and all inequalities become equalities.

Therefore $u_{i,i-1} \geq 1$ and $dr=\sum_{i=2}^dr
u_{i,i-1}+u_{1,dr} \geq dr-1+u_{1,dr}$, which implies that
$u_{1,dr}=1$. As $\ffi$ decreases from top to bottom and bottom to
left, we have that $u_{i,i-1}=1$. Since $\ffi$ is homogeneous we can
conclude that all the entries in $\ffi$ have degree 1.

We prove b).

(i) $\Rightarrow$ (ii) If $\EE$ has a  linear minimal resolution $$0 \ra \OO_{\PP^n}(t-1)^s
\srel{\ffi}{\ra} \OO_{\PP^n}(t)^s \ra  \EE \ra 0 ,$$ then $\ffi$ is a reduced
matrix factorization with $\det (\ffi)=f^r$, $r=\rg(\EE)$ and the
degree of $\det (\ffi)$ is $s$. Hence $s= \deg(X)r$ and $\EE$ is Ulrich.
Moreover as $\EE$ is normalized, this implies that $t=0$ and $\EE$ is generated by its global sections.

(ii) $\Rightarrow$ (i)
Let $$0 \lra \FF=\oplus_{j=1}^{dr} \OO_{\PP^n}(b_j)
\srel{\ffi}{\lra} \GG=\oplus_{i=1}^{dr} \OO_{\PP^n}(a_i) \lra \EE \lra
0$$ be the minimal free resolution of $\EE$ as $\OO_{\PP^n}$-module, with $a_1
\geq \dots \geq a_{dr}$ and $b_1 \geq \dots \geq b_{dr}$,
$(\varphi,\psi)$  a reduced matrix factorization of $f$ (note that
by \ref{lemEUlrich} $\EE$ has no free direct summand). As $H^0(\EE(-1))=0$, we
have $H^0(\FF(-1)) \cong H^0(\GG(-1))$. Therefore we have an
equality $\{b_j \mid b_j \geq 1\}=\{a_i \mid a_i \geq 1\}$, and in
particular $a_1=b_1$ if these sets are non-empty. But then, as
$\ffi$ is reduced, its last column would be zero which is
impossible because $\ffi$ is a matrix factorization. Hence, $a_i
\leq 0$ and $b_i < 0$ for all $i=1, \dots ,m$ and as $H^0(\EE)
\neq 0$, $a_1 =0$.

Now, $\EE$ is generated by global sections of and only if $a_{dr}=0$,
which is now equivalent to saying that all $a_i$ are zero. But the
the degree of the determinant of $\ffi$ must be $dr$, so this is
equivalent to saying that $\EE$ has a linear resolution.

(iii) $\Rightarrow$ (i) We assume now that $\EE$ is  normalized with $h^0(\EE)=dr$. Then
by Corollary \ref{boundh0} $\EE$ is Ulrich. To see that is has a
linear minimal free resolution we proceed by induction on
$r=\rg(\EE)$. For $r=1$, $\EE$ is not extension of lower rank
Ulrich sheaves, so we conclude by a).

We assume now that $r>1$. If $\EE$ is not extension of lower rank
Ulrich sheaves, then by a) we are done. If $\EE$ is an
extension of Ulrich sheaves $\EE_i$ of lower rank $r_i$, $i=1,2$,
$$0 \lra \EE_1 \lra \EE \lra \EE_2 \lra 0,$$
then $h^0(\EE_i)= 3r_i$ because of the bound given in Corollary
\ref{boundh0}. In particular, $\EE_i$ is Ulrich and normalized and
by induction hypothesis we  obtain that $\EE_i$ has a linear
minimal resolution. We then apply horseshoe lemma in \cite[2.2.8]{Weibel}
 to obtain a linear resolution for $\EE$.

%Moreover, dualizing the exact
%sequence
%$$0 \lra E(-d) \lra R_X(-1)^{dr} \lra E^{\sigma} \lra 0$$
%obtained from the resolution in (iii), we have the exact sequence
%$$0 \lra E^{\sigma\vee} \lra R_X(1)^{dr} \lra E^{\vee}(d) \lra 0
%.$$ Hence, $E^{\vee}(-d-1)$ is normalized.
\end{proof}

The existence of Ulrich sheaves on a projective scheme $X$ is not known in general. Below we
present some known examples of Ulrich sheaves.

\begin{ex}\rm a) The
existence of an Ulrich sheaf of rank $1$ on a hypersurface defined by a form $f$ is related to the
possibility of writing $f$ as the determinant of a $d \times
d$-matrix of linear forms. For example, on a smooth cubic surface
the existence of twisted cubic curves in $X$ allows writing the
equation of $f$ as a linear determinant (see \cite{Beauville}
Corollary 6.4).

b) The existence orientable rank 2 Ulrich bundles on a hypersurface defined by a form $f$ is related to
the possibility of writing $f$ as the Pfaffian of a linear
skew-symmetric matrix, and by Serre's correspondence, to the
existence of certain arithmetically Gorenstein subschemes of
codimension 2 of $X$. On a smooth cubic surface, the Serre
correspondence applied to a set of 5 general  points in $X$ proves
the existence of a rank 2 Ulrich bundle with a linear resolution
(\cite[Proposition 7.6]{Beauville}).

c) The existence of Ulrich sheaves on special varieties has been studied by many authors. It is known for example that there exists at least one Ulrich MCM module on complete intersections (see \cite{HerzogUlrichBackelin} and the references therein). In \cite{EisSchreyer}, Eisenbud and Schreyer proved the existence of Ulrich bundles on any algebraic curve, on Veronese varieties and 
%they also prove the existence of 
of rank two Ulrich sheaves 
on Del Pezzo surfaces. Although rank one Ulrich sheaves are rare in general, it is not difficult to prove that there exists a rank one Ulrich sheaf on any ACM rational surface $S$ in $\PP^4$ (in this case $S$ is either a cubic scroll, a Del Pezzo surface, a Castelnuovo surface or a Bordiga surface, and the existence of a rational quartic curve on them leads to the existence of a rank one Ulrich sheaf).
%It is known that there exists at least one Ulrich MCM module over certain local rings such as hypersurface rings and complete intersection rings  (see \cite{HerzogUlrichBackelin} and the references therein).
\end{ex}

%\begin{prop}\label{dualUlrich} Let $\EE$ be an ACM  sheaf on  a hypersurface $X
%\subseteq \PP^n$. Then
%\end{prop}
%
%\begin{proof} If $\EE$ is Ulrich then by \ref{resolUlrich}, its
%minimal free $\OO_X$-resolution is 2-periodic and we can split it
%into short exact sequences:
%
%DIAGRAM $$\cdots \lra \OO_X(-d)^{dr} \srel{\psi}{\lra}
%\OO_X(-1)^{dr} \srel{\varphi}\lra \OO_X^{dr} \lra \EE \lra 0.$$
%
%Therefore $\EE=\ker(\varphi)$. As $\mathcal{E}xt^1(\EE,\OO_X)=0$
%because $\EE$ is locally CM, taking $\mathcal{H}om(\,\, ,\OO_X)$
%in the the long exact sequence is exact. In particular, we obtain
%that $\EE^{\vee}=\coker(\varphi^{\vee})$ and has  minimal free
%resolution starting as
%$$\cdots \lra \OO_X^{dr} \srel{\varphi^{\vee}}\lra \OO_X(1)^{dr} \lra \EE^{\vee} \lra 0.$$
%Therefore it is an Ulrich sheaf.
%
% RATHER WITH MODULES TO ENSURE MINIMAL REsOLUTION???
%
%Conversely if $\EE^{\vee}$ is Ulrich, we just need to change the
%roles of $\EE$ and $\EE^{\vee}$ in the previous argument as $\EE$
%is reflexive.
%
%\end{proof}

\section{Ulrich vector bundles on the cubic surface}
 From now on, $X$ will be a non-singular cubic surface in $\PP^3$ defined by a degree $3$ homogeneous polynomial $f \in R=k[x_0,x_1,x_2,x_3]$
 and $R_X$ will denote the ring $R/(f)$.
 Then $p_a(X)=0$ and $\omega_X=\OO_X(-1)$. Moreover, as $X$ is regular,
any ACM sheaf on $X$ is locally free.
Let $\EE$ be a vector bundle of rank $r$ on $X$. Let  $c_1(\EE)$ and $c_2(\EE)$ denote its Chern classes,  and $\deg(\EE)=\deg c_1(\EE)$  its degree. Then the Riemann-Roch theorem
says that
$$\chi(\EE)=r+\frac{c_1(\EE).H}{2}+ \frac{c_1(\EE)^2-2c_2(\EE)}{2} \,.$$
$\EE$ is called \emph{orientable} if $\det(\EE)$ is isomorphic to
$\OO_X(l)$ for some $l\in \ZZ$, or in other words, $c_1(\EE)=lH$.

\begin{rk}\label{bounddeg}\rm Let $\EE$ be  an ACM bundle of rank $r$ on  the cubic surface $X$, and let $H$ be a general hyperplane section. Then from the exact sequence
$$0 \lra \EE(-1) \lra \EE \lra \EE_H \lra 0$$ we find $h^0(\EE) \geq h^0(\EE_H)$ and $\deg \EE=\deg \EE_H$. By Riemann-Roch on the elliptic cubic curve
we have $h^0(\EE_H)\geq \deg \EE_H$ and so $\deg\EE \leq
h^0(\EE)$.

On the other hand, we know by Corollary \ref{boundh0} that if
$\EE$ is ACM and normalized, $h^0(\EE) \leq 3r$. Therefore, if
$\EE$ is a normalized ACM bundle on $X$, then  $\deg(\EE)\leq 3r$
and equality implies that $\EE$ is an Ulrich bundle with
$h^0(\EE)=3r$.
\end{rk}

\begin{lema}\label{chernUlrich} Let $\EE$ be a normalized ACM bundle of rank $r$ on $X$.
Then the following are equivalent
\begin{enumerate}
\item $\EE$ has a linear minimal  free resolution
\item $\EE$ is Ulrich and generated by its global sections
\item $h^0(\EE)=3r$
\item $\deg(\EE)=3r$
\end{enumerate}
%Moreover, if these conditions hold, then $c_2(\EE) \leq
%\frac{1}{2}(3r^2-r)$ and  we have equality if $\EE$ is orientable.
\end{lema}

\begin{proof}
(i), (ii) and (iii) are equivalent by  Lemma \ref{resolUlrich}.

(i) $\Rightarrow$ (iv): From the exact sequence $$0 \lra
\OO_{\PP^3}(-1)^{3r} \lra \OO_{\PP^3}^{3r} \lra \EE \lra 0$$ we
obtain that $\chi(\EE)=3r$ and $\chi(\EE(1))=9r$. For a general
hyperplane $H$, the hyperplane section of $X$ is an integral
elliptic curve and we have the exact sequence
$$0 \lra \EE \lra \EE(1) \lra \EE_H(1) \lra 0.$$
Therefore $\chi(\EE(1))-\chi(\EE)=\chi(\EE_H(1))$. By Riemann-Roch
theorem on a elliptic curve, this last term is equal to
$\deg(\EE_H(1))$. As $c_1(\EE_H(1))=c_1(\EE_H)+rH$, we have that
$\deg(\EE_H(1))=\deg(\EE)+3r$ and this implies that
$\chi(\EE(1))-\chi(\EE)=\deg(\EE)+3r$. As the term in the left is
equal to $6r$, we obtain $\deg(\EE)=3r$ as desired.

(iv)  $\Rightarrow$  (ii) by Remark \ref{bounddeg}. Therefore all
conditions are equivalent.
\end{proof}

We recall the following result that was proved by the first author in \cite{Cmrc}.
\begin{prop}[see \cite{Cmrc}]\label{resolpoints}
Let $Z$  be a set of $n=\frac{1}{2}(3r^2-r)$ general
points on $X$, $r \geq 2$. Then the minimal free resolution over $R$ of the
saturated ideal of $Z$ in $X$  is
$$0 \ra R(-r-3)^{r-1} \ra R(-r-1)^{3r}  \ra R(-r)^{2r+1}   \ra I_{Z,X}
 \ra 0.$$
% Then the minimal free resolution of the
% saturated ideal of $Z$ in $R$ is
% $$0 \ra R(-r-3)^{r-1} \ra R(-r-1)^{3r}  \ra R(-r)^{2r+1} \oplus R(-3)  \ra I_{Z}
%  \ra 0.$$
\end{prop}

In the next theorem  we prove the existence of Ulrich bundles of
any rank.
\begin{teo}\label{generalpoints}
a) If $ \EE$ is a normalized orientable Ulrich bundle of rank $r \geq 2$ generated by global sections, then there is an exact sequence
\begin{equation}\label{extTeo}
0 \lra \OO_X^{r-1} \lra \EE \lra \JJ_{Z,X}(r) \lra 0
\end{equation}
where $Z$ is a zero-scheme of degree $n=\frac{1}{2}(3r^2-r)$, and $h^0(\JJ_{Z,X}(r-1))=0$.

b) Conversely, if $r\geq 2$ and $Z$ is a sufficiently general set of $n=\frac{1}{2}(3r^2-r)$ points on $X$, then there exists an extension of $\JJ_{Z,X}(r)$ by $\OO_X^{r-1}$ as above, where $\EE$ is a normalized orientable Ulrich bundle of rank $r$ generated by global sections.

%  and  the ideal sheaf of $Z$ in $X$ has an $\mathcal{N}$- type
%  resolution
% $$0 \lra \OO_X^{r-1} \lra \UU_r \lra \JJ_{Z,X}(r) \lra 0 $$
% where $\UU_r$ is a normalized orientable Ulrich bundle of rank $r$
% with $c_1(\UU_r)=rH$ and $c_2=n$.
\end{teo}

\begin{proof}
a) As $\EE$ is generated by
global sections, we can take $r-1$ general sections of $\EE$ so
that the quotient of $\EE$ by $\OO_X^{r-1}$ is torsion free.
Therefore we have an exact sequence
$$0 \lra \OO_X^{r-1}
\lra \EE \lra \JJ_Z(D) \lra 0$$ where $D=c_1(\EE)$ is a certain
divisor on $X$ and $Z$ is a zero-scheme of degree equal to
$c_2(\EE)$.
%The divisor $D$ is effective and
By Lemma \ref{chernUlrich}(iv) $D$ has degree
$3r$. As $\EE$ is orientable we have $c_1(\EE)=rH$, and the Riemann-Roch theorem applied to $\EE$ gives $c_2(\EE)
=\frac{1}{2}(3r^2-r)$. Note that $h^0(\JJ_{Z,X}(r-1))=0$ since $\EE$ is normalized.
% \begin{equation}\label{eqRR}
% c_2(\EE)=-2r+\frac{1}{2}D^2+\frac{1}{2}3r \,.
% \end{equation}
%
% We bound $D^2$ as follows. As $D$ is a curve on $X$ of degree
% $3r$, its genus is bounded above by the genus of a complete $C$
% intersection of degree $3r$, $C \sim rH$ (see \cite{AG} Ex
% III.4.7). In particular, $2g_D -2 \leq 2g_C-2$ and by the
% adjunction formula we obtain $D^2-\deg D \leq C^2 -\deg C$.
% Therefore, $D^2 \leq C^2=3r^2$ and from equation (\ref{eqRR}) we
% get
% $$c_2(\EE) \leq \frac{1}{2}(3r^2-r) \,.$$
%
% Moreover, as $\EE$ is orientable, we have $c_1(\EE)=rH$ and by Riemann-Roch we have $c_2(\EE)
% =\frac{1}{2}(3r^2-r)$.

b) Let $Z$ be a set of $n$ points, sufficiently general so that the  minimal free resolution of $\JJ_{Z,X}$ is given by Proposition \ref{resolpoints}.
By a generalization of the well-known Serre correspondence we obtain an extension of $\JJ_{Z,X}(r)$ by $\OO_X^{r-1}$ in the following way.

If $R_X$, $R_Z$ denote the homogeneous coordinate rings of $X$ and $Z$, and $I_{Z,X}$ the saturated ideal of $Z$ in $X$, then from the exact sequence
$$0 \lra I_{Z,X} \lra R_X \lra R_Z \lra 0$$
we obtain $Ext^1(I_{Z,X},R_X(-1)) \cong Ext^2(R_Z, R_X(-1))$. As the canonical module of $X$ is isomorphic to $R_X(-1)$, we have that $Ext^2(R_Z, R_X(-1))$ is isomorphic to the canonical module  $K_Z$ of $Z$, because 2 is precisely the codimension of $Z$ in $X$.

On the other hand, a minimal free resolution for the canonical module $K_Z$ can be obtained by applying $Hom_R(\cdot,K_{\PP^3})$ to  the  minimal free resolution of $I_Z$ given in Proposition \ref{resolpoints} (see \cite[1.2.4]{M}). Therefore there is a minimal resolution of $K_Z$ as follows:
$$ \dots \lra R(r-3)^{3r} \lra R(r-1)^{r-1} \lra K_Z \lra 0 \,.$$
 and $K_Z$ is generated in degree $1-r$ by elements $\nu_1, \dots, \nu_{r-1}$.
These generators provide an extension
\begin{equation}\label{ext}
0 \lra R_X^{r-1} \lra E \lra I_{Z,X}(r) \lra 0
\end{equation}
via the isomorphism $K_Z \cong Ext^1(I_{Z,X},R_X(-1))$. To prove that this module $E$ is a maximal Cohen-Macaulay $R_X$-module we just need to prove that $Ext^1(E,K_X)=0$. This follows by applying $Hom_{R_X}(\cdot,K_{X})$ to the exact sequence \eqref{ext}. Indeed, this leads to an exact sequence
$$  Hom(R_X^{r-1},K_X) \cong R_X(-1)^{r-1} \ra Ext^1(I_{Z,X}(r),K_X)\cong K_Z(-r) \ra Ext^1(E,K_X) \ra 0$$
where the first morphism is an epimorphism because it is defined by the generators $\nu_1, \dots, \nu_{r-1}$.

Now observe that $\EE:=\widetilde{E}$ is normalized. Indeed, it can be deduced from the exact sequence \eqref{ext} that $h^0{\EE}=3r$ and that $h^0(\EE(-1))=0$ ($h^0(\JJ_{Z,X}(r))$ and $h^0(\JJ_{Z,X}(r-1))$ can be computed from the resolution given in Proposition \ref{resolpoints}). Therefore, by Lemma \ref{chernUlrich} we obtain that $\EE$ is Ulrich and generated by global sections.
\end{proof}

\begin{cor}
On a nonsingular cubic surface $X \subseteq \PP^3$, for every $r \geq 2$, there exist normalized orientable Ulrich bundles $\EE$ of rank $r$ generated by global sections with $c_1(\EE)=rH$ and $c_2(\EE)=\frac{1}{2}(3r^2-r)$
\end{cor}

\section{Stability of general Ulrich bundles}
Our goal in this section is to prove that the Ulrich bundles constructed
in the previous section are stable. Following the terminology of \cite{HL}, we recall that a vector bundle $\EE$ on a nonsingular projective variety $X$ is \textit{semistable} if for every coherent subsheaf $\FF $ of $\EE$ we have the inequality
$$P(\FF)/\rg(\FF) \leq P(\EE)/\rg(\EE),$$
where $P(\FF)$ and $P(\EE)$ are the Hilbert polynomials of the sheaves. It is \textit{stable} if one always has strict inequality above. With these definitions, one knows (cf.  \cite[4.3.4]{HL}) that there is a projective coarse moduli scheme $M^{ss}(P)$ whose closed points are in one-to-one correspondence with certain equivalence classes of semistable sheaves, and there is an open subscheme $M^s(P)$ whose points correspond to the isomorphism classes of stable vector bundles.

There is another definition, more adapted to calculations, using
the \textit{slope} of $\EE$, which is defined as
$\deg(c_1(\EE))/\rg(\EE)$, where $c_1(\EE)$ is the first Chern
class. We say that $\EE$ is $\mu$-\textit{semistable} if for every
subsheaf $\FF$ of $\EE$ with $0 < \rg \FF < \rg \EE $, $slope
(\FF) \leq slope (\EE)$. We say $\EE$ is $\mu$-\textit{stable} is
strict inequality $<$ always holds. The two definitions are
related as follows
$$\mu-stable  \quad \Rightarrow \quad stable \quad \Rightarrow \quad semistable \quad \Rightarrow \quad \mu-semistable$$
% To this end we need to make a
To begin with, we need a
summary of Atiyah's classification of vector bundles on an
elliptic curve.

\begin{rk}\label{Atiyah}\rm \textbf{(Vector bundles on an elliptic
curve.)}
\\
\noindent Atiyah proved in \cite{Atiyah} that the following statements hold
for vector bundles on a nonsingular elliptic curve $Y$.
\begin{itemize}
\item For every $r \geq 1$, $d \in \ZZ$, there is a 1-dimensional
family of indecomposable bundles $\EE_{r,d}$ of rank $r$ and
degree $d$ parameterized by the points of $Y$.
\item $h^0(\EE_{r,d})=
  \begin{cases}
    d & \text{ if } d>0, \\
    0 \text{ or } 1 & \text{ if } d=0\\
    0 & \text{ if } d<0
  \end{cases}.
$
\item $\EE_{r,d}$ is semistable for all $r,d$ and it is stable if
$(r,d)=1$.
\item Every vector bundle on $Y$ is a direct sum of indecomposable
bundles $\EE_{r,d}$.
\end{itemize}
Now let $Y \subset \PP^2$ be a nonsingular cubic curve. From the
above it is easy to see that the only normalized Ulrich bundles on
$Y$ are $\EE_{r,d}$ with $d=3r$ in the case for which $h^0(
\EE_{r,d}(-1))=0$. These correspond to an open subset of the curve
$Y$. They are semistable but not stable and satisfy $h^0(
\EE_{r,d})=3r.$ (Note that on a curve the two definitions of stability and semistability coincide).
\end{rk}

We come back to the nonsingular cubic surface $X$.

%We suppose
%given a bundle $\EE$ satisfying
%\begin{center} ($\star$) $\EE$ is ACM of rank $r$, and admits an exact
%sequence $$0 \lra \OO_X^{r-1} \lra \EE \lra \JJ_Z(r) \lra 0$$
%where $Z$ is a set of $n=\frac{1}{2}(3r^2-r)$ points and
%$h^0(\JJ_Z(r-1))=0.$
%\end{center}
%Then $\EE$ is normalized and has degree $3r$, hence by
%\ref{chernUlrich} it is an Ulrich bundle generated by global
%sections and satisfies $h^0(\EE)=3r$. Moreover it has slope 3.

\begin{prop}\label{propstable} Any bundle $\EE$ satisfying \ref{generalpoints} a) is $\mu$-semistable.
\end{prop}
\begin{proof} Since $\EE$ is normalized and ACM with $h^0(\EE)=3r$ it
follows that the general hyperplane section $\EE_H$ is also
normalized of degree $3r$ and $h^0(\EE_H)=3r$. By Atiyah's
classification $\EE_H=\oplus_i \EE_{r_i,d_i}$. Since
$h^0(\EE_{r_i,d_i}(-1))=0$ it follows that $d_i \leq 3r_i$. Summing
up, we must have equality for each $i$. Since the bundles
$\EE_{r,d}$ on the elliptic curve are all semistable, we find that
$\EE_H$ is semistable of slope 3.

Now if $\FF$ is a coherent subsheaf of $\EE$, $\FF_H$ is a
coherent subsheaf of $\EE_H$, so $\FF_H$ has slope $\leq 3$. Hence
$\FF$ has slope $\leq 3$ and $\EE$ is $\mu$-semistable.
\end{proof}

\begin{teo}\label{stable} If $\EE$ is a vector bundle on $X$ satisfying \ref{generalpoints} a) and $Z$ is sufficiently general,  $\EE$ is $\mu$-stable.
\end{teo}

\noindent\textit{Proof} Let $\FF \subset \EE$ be a coherent
subsheaf. Since $\EE$ is $\mu$-semistable, $slope(\FF) \leq
slope(\EE)=3$. We need only to eliminate the case $slope(\FF)=3$.

By pulling-back torsion if necessary, we may assume that $\EE/\FF$
is torsion-free, in which case $\FF$ is locally free. So we only
need to show the nonexistence of a semistable locally free proper
subsheaf $\FF$ of slope 3.

From the inclusion $\FF \subset \EE$ we find $Hom(\FF,\EE)\neq 0$,
so $H^0(\FF^{\vee} \otimes \EE) \neq 0$. Tensoring the sequence
\eqref{extTeo} of \ref{generalpoints}a) with $\FF^{\vee}$ we find
$$0 \lra (\FF^{\vee})^{r-1} \lra \EE \otimes \FF^{\vee} \lra \JJ_Z(r) \otimes \FF^{\vee} \lra 0.$$
Since $\FF^{\vee}$ is $\mu$-semistable of negative degree, it has no
sections. So we will achieve a contradiction by showing that
$h^0(\JJ_Z(r) \otimes \FF^{\vee})=0$ for $Z$ sufficiently general.
At this point we need a lemma.

\begin{lema} Let $\GG$ be a $\mu$-semistable vector bundle on $X$ of
rank $s$ and degree $3sp$ for some $p>0$. Suppose also that
$\GG_H$ is semistable for a general plane section. Then $h^0(\GG)
\leq \frac{3}{2}sp(p+1)+s.$
\end{lema}
\begin{proof} From the exact sequence
$$0 \lra \GG(-1) \lra \GG \lra \GG_H \lra 0$$
we have $h^0(\GG) \leq h^0(\GG(-1)) + h^0(\GG_H)$. Applying the
same to $\GG(-1), \GG(-2) \dots$ and summing we find $$h^0(\GG)
\leq  \sum_{i=0}^{\infty} h^0(\GG(-i)).$$

Note that $\deg(\GG(-i))=3s(p-i)$. Since $\GG_H(-i)$ is
semistable, all its indecomposable summands have the same slope.
Then using Atiyah's results we obtain $$\begin{array}{ll}
    h^0(\GG(-i))=3s(p-i) & \text{ for } 3s(p-i)>0, \\
    h^0(\GG(-i)) \leq \rg=s & \text{ for } 3s(p-i)=0,\\
     h^0(\GG(-i)) =0 & \text{ for } 3s(p-i)<0,
\end{array}$$
Hence $h^0(\GG) \leq 3s(1+2+ \dots +p)+s= \frac{3}{2}sp(p+1)+s$.
\end{proof}

\begin{proofteo} We apply the lemma to $\FF^{\vee}(r)$ which has
rank $s < r$ and degree $-3s+3rs=3s(r-1)$. Hence
$h^0(\FF^{\vee}(r)) \leq \frac{3}{2}sr(r-1)+s$.

If $H^0(\FF^{\vee}(r))=0$, we are done.

If $H^0(\FF^{\vee}(r))\neq 0$, suppose first that  $H^0(\FF^{\vee}(r))$ generates a subsheaf of
rank $s$ of $\FF^{\vee}(r)$. Then each general point $P$ will
impose $s$ conditions on a section of $\FF^{\vee}(r)$ to lie in
$\JJ_P \otimes \FF^{\vee}(r)$. Since $Z$ consists of
$n=\frac{1}{2}(3r^2-r)$ general points, we have
$\frac{1}{2}s(3r^2-r)$ potential conditions. This number is
greater than $h^0(\FF^{\vee}(r))$, so we conclude that $h^0(\JJ_Z
\otimes \FF^{\vee}(r))=0$ as required.

Suppose on the other hand that the sections of $\FF^{\vee}(r)$
generate a subsheaf $\GG$ of rank $t <s$. Then each general point
$P$ will impose only $t$ conditions on a section  to lie in $\JJ_P
\otimes \FF^{\vee}(r)$. But $\GG$ and $\GG_H$ will also be
$\mu$-semistable of slope $\leq 3$, so the arguments of the lemma will
give a similar bound with $t$ in place of $s$, and we conclude as
before.
\end{proofteo}

Notice that usually the existence of ACM bundles of high rank is only proved using extensions of lower rank ACM bundles.
In the following result we show that the bundles constructed in Theorem \ref{generalpoints} (or more generally, $\mu$-stable Ulrich bundles on $X$) are not extension of ACM bundles of lower rank.
\begin{prop}\label{indec} Let $\EE$ be a $\mu$-stable Ulrich bundle of rank $r$ with $c_1=rH$. Then $\EE$ is normalized and it is not an extension of lower rank ACM bundles. In particular, the bundles constructed in \ref{generalpoints} b) are indecomposable and are not extensions of lower rank ACM bundles.
\end{prop}

\begin{proof}
We first prove that $\EE$ is normalized.

As $\EE$ is Ulrich and $\mu$-stable, it is easy to see that it cannot be an extension of Ulrich bundles of lower rank. Therefore by Proposition  \ref{resolUlrich}(1), $\EE$ has a linear resolution
$$0 \lra \OO_{\PP^n}(a-1)^{3r} \srel{\varphi}\lra \OO_{\PP^n}^{3r}(a) \lra \EE \lra 0.$$
Then $\chi(\EE)=3r\left({\binom{a+3}{3}}-{\binom{a+2}{3}}\right)$ and $\chi(\EE(-1))=3r\left({\binom{a+2}{3}}-{\binom{a+1}{3}}\right)$. The difference between these two values is $\chi(\EE)-\chi(\EE(-1))=3r(a+1)$.
On the other hand, the exact sequence
$$0 \lra \EE(-1) \lra \EE \lra \EE_H \lra 0$$
tells us that $\chi(\EE)-\chi(\EE(-1))=\chi(\EE_H)$ and this last term is equal to $\deg\EE$ by Riemann-Roch on the elliptic curve. We conclude that $a=0$ and $\EE$ is normalized.

Now suppose there exists an extension $$0\lra \EE_1 \lra
\EE \lra \EE_2 \lra 0$$ with $\EE_i$ ACM bundle of rank $r_i$,
$r=r_1+r_2$. As $\EE$ is $\mu$-stable, the slope of $\EE_1$ is $<3$ so
$\deg(\EE_1)<3r_1$. On the other hand, by Remark \ref{bounddeg} we
have $\deg(\EE_2) \leq h^0(\EE_2)$ and by Corollary \ref{boundh0},
$h^0(\EE_2) \leq 3r_2$ (this last bound still hold in case
$h^0(\EE_2)=0$). Therefore $3r=3r_1+3r_2
> \deg(\EE_1)+\deg(\EE_2)=\deg(\EE)$ which is precisely $3r$ and
we obtain a contradiction.
\end{proof}

\begin{lema} With $\EE$ as in \ref{generalpoints} a) we have $\chi(\EE \otimes
\EE^{\vee})=-r^2$.
\end{lema}
\begin{proof} Use additivity of $\chi$ on the sequences
\begin{equation*}0 \lra (\EE^{\vee})^{r-1} \lra \EE \otimes \EE^{\vee} \lra
\JJ_Z(r) \otimes \EE^{\vee} \lra 0\end{equation*}
\begin{equation*}0 \lra \OO_X(-r) \lra  \EE^{\vee} \lra
\OO_X^{r-1}  \lra \EE xt^1(\JJ_Z(r), \OO_X) \lra 0\end{equation*}
\begin{equation*}0 \lra \JJ_Z(r) \otimes\EE^{\vee} \lra \EE^{\vee}(r)  \lra
 \EE^{\vee}(r) \otimes \OO_Z \lra 0\end{equation*}
\begin{equation*}0 \lra \OO_X \lra  \EE^{\vee}(r) \lra
\OO_X^{r-1}(r)  \lra \EE xt^1(\JJ_Z(r), \OO_X(r)) \lra
0\end{equation*} and the fact that the last sheaf in the second
and fourth sequences has length $n$ and the last sheaf in the
third sequence has length $rn$. The calculations are left to the
reader.
\end{proof}

%\begin{cor} An Ulrich bundle $\UU_r$ corresponding to a set of general
%points
Summing up we obtain the main result of this paper.
\begin{teo}\label{mainteo} Let $X$ be a nonsingular cubic surface in $\PP^3$ over an algebraically closed field $k$. Let $r\geq 2$. The normalized orientable Ulrich bundles of rank $r$ generated by global sections are all $\mu$-semistable.
% The coarse moduli space is irreducible and has a non-empty smooth open subset of dimension $r^2+1$ corresponding to stable vector bundles.
Among these, the $\mu$-stable Ulrich bundles correspond to a non-empty open subset of dimension $r^2 +1$ of the moduli space $M^s(r;c_1,c_2)$ of stable vector bundles on $X$ with Chern classes $c_1=rH$ and $c_2=\frac{1}{2}(3r^2-r)$.
\end{teo}

\begin{proof} The first statement follows from \ref{propstable}. If $\EE$ is $\mu$-stable then it is stable and $h^0(\EE \otimes
\EE^{\vee})=1$. By duality $h^2(\EE \otimes \EE^{\vee})=h^0(\EE
\otimes \EE^{\vee}(-1))=0$. Hence from the previous lemma we find
$h^1(\EE \otimes \EE^{\vee})=r^2+1$.

Since $h^2(\EE \otimes \EE^{\vee})=0$, there are no obstructions, so at the point corresponding to $\EE$,
the moduli space is smooth of dimension $h^1(\EE \otimes
\EE^{\vee})$ (cf. \cite[4.5.2]{HL}).

% In the following corollary we prove that the bundles considered so far form an open subset of the moduli space of stable vector bundles.
%
% \begin{cor}
% The moduli space $M^s_X(r;c_1,c_2)$ of rank $r$ stable vector bundles on $X$ with Chern classes $c_1 =rH$, $c_2=\frac{1}{2}(3r^2-r)$ has an irreducible component of dimension $r^2+1$ containing a smooth open subset that corresponds to stable Ulrich bundles.
% \end{cor}
% \begin{proof}
It remains to show that the $\mu$-stable Ulrich bundles with the given Chern classes form an open subset of the moduli space. This open set will be nonempty by \eqref{generalpoints} and \eqref{stable}.

In any flat family of vector bundles, the condition of being ACM is an open condition by the semicontinuity theorem \cite[III, 12.8]{AG}. The condition of being Ulrich is not necessarily open, but the condition of being $\mu$-stable is open and here we prove that $\mu$-stable ACM bundles with the given Chern classes are Ulrich.

Indeed, let $\EE$ be an ACM $\mu$-stable bundle on $X$  with $c_1 =rH$, $c_2=\frac{1}{2}(3r^2-r)$. As $\EE$ is $\mu$-stable of slope 3, $\EE(-1)$ cannot have sections.  $\EE^{\vee}$ is also $\mu$-stable of slope $-3$, so we have $h^0(\EE^{\vee}(-1))=0$. By Serre duality we obtain $h^2(\EE)=0$ and then by Riemann-Roch theorem on $X$ it follows that $h^0(\EE)=3r$.
Thus $\EE$ is normalized with $h^0(\EE)=3r$ and therefore it is Ulrich by Lemma \ref{chernUlrich}.
\end{proof}

\section{Filtrations of general ACM bundles}

Although we have seen in Proposition \ref{indec} the the bundles
$\EE$ constructed in Theorem \ref{generalpoints} are not extension
of lower rank ACM bundles, we shall prove that for a suitable
dissoci\'e sheaf $\LL$, the bundle $\EE \oplus \LL$ is indeed an
extension of ACM bundles of rank 1.

Our main technique in this section is Gorenstein liaison theory
or, more precisely, strict $G$-links (see definition \ref{link}
below). In this case the links are performed by divisor of type
$mH_Y-K_Y$ on an ACM scheme $Y$ that is a divisor on a scheme $X
\subseteq \PP^n$. To put this in a more general context, recall
that if $Y$ is an ACM scheme satisfying $G_0$, we can define the
{\em anticanonical divisor} $M = M_Y$, given by an embedding of
$\omega_Y$ as a fractional ideal in the sheaf of total quotient
rings ${\mathcal K}_Y$, even if $Y$ does not have a well-defined
canonical divisor \cite[2.7]{Har_divisors_biliaison}. Recall also
that if $Y$ is an ACM scheme in ${\mathbb P}^n$ satisfying $G_0$,
and if $G$ is an effective divisor on $Y$, linearly equivalent to
$mH+M_Y$, for some $m \in {\mathbb Z}$, then $G$ is AG in
${\mathbb P}^n$ \cite[3.4]{Har_divisors_biliaison},
\cite[5.2,5.4]{KMMNP}.

\begin{defi}\label{link} Let $Y \subset \PP^n$ be an arithmetically Cohen-Macaulay scheme that is Gorenstein in codimension 0.
We say that two equidimensional subschemes $Z_1, Z_2$ of
codimension 1 of $Y$ without embedded components  are
\textit{strictly} $G$-\textit{linked} by a subscheme $G \subset Y$
%if $Z_1+Z_2=G$ and
if $G$ contains $Z_1$ and $Z_2$, $\JJ_{Z_i} \cong \mathcal{H}om
(\OO_{Z_j}, \OO)$ for  $i,j=1,2$, $i \neq j$, and there is an
$m\in \ZZ$ such that $G$ is linearly equivalent to $mH_Y+M_Y$
where $M_Y$ is the anticanonical divisor
\cite{Har_divisors_biliaison}.
\end{defi}
In our case, $Y$ will be an ACM curve contained in the nonsingular
cubic surface $X$ and therefore it will have a well-defined
canonical divisor $K_Y$ so that $M_Y=-K_Y$. Divisors linearly
equivalent to $mH-K_Y$ are AG schemes \cite[3.3]{Har_divisors_biliaison}.

To associate ACM bundles to sets of points we work with $\NN$-type
resolutions (or so-called Bourbaki sequences, cf. \cite{BHU}). We
recall the definition here.

\begin{defi} Let $X \subset \PP^n$ be an equidimensional scheme and
$Z \subset X$ be  a codimension 2 subscheme without embedded
components. An $\NN$-\textit{type resolution} of $Z$ is an exact
sequence
\[ 0 \rightarrow {\mathcal L} \rightarrow {\mathcal N}
\rightarrow {\JJ}_{Z,X} \rightarrow 0
\]
with ${\mathcal L}$ dissoci\'e and ${\mathcal N}$ a coherent sheaf
satisfying $H_*^1({\mathcal N}^{\vee}) = 0$ and ${\mathcal
E}xt^1({\mathcal N},{\mathcal O}_X) = 0.$
\end{defi}
When $X$ satisfies Serre's condition $S_2$ and $H_*^1({\mathcal
O}_X) = 0$, then an ${\mathcal N}$-type resolution of $Z$ exists
(\cite[2.12]{Hart_LR}). An $\NN$-type resolution is not unique but
it is well known that any two $\NN$-type resolutions of the same
subscheme are \textit{stably equivalent} (see
\cite[1.10]{Hart_LR}). In other words, if $\NN$ and $\NN'$ are two
sheaves appearing in the middle of an $\NN$-type resolution of a
subscheme $Z$, then there exist dissoci\'e sheaves $\LL_1$,
$\LL_2$ and an integer $a$ such that
$$\NN \oplus \LL_1 \cong \NN'(a) \oplus \LL_2.$$ See \ref{ex2} and
\ref{ex1} below for examples of $\NN$-type resolutions.

When $X$ is an AG scheme and we have an $\NN$-type resolution as
above, then $\NN$ is an ACM sheaf if and only if $Z$ is an ACM
scheme. In particular, if $Z$ is a 0-dimensional scheme, then any
sheaf appearing in an $\NN$-type resolution of $Z$ will be ACM. We
already saw examples of  $\NN$-type resolutions in Theorem
\ref{generalpoints}.

We study when an AG scheme $G$ of codimension 2 in an AG scheme $X
\subset \PP^n$ occurs as a divisor $mH_Y-K_Y$ on some ACM divisor
$Y \subset X$.
\begin{prop}\label{Gorenstein} Let $X$ be an AG
scheme with $\omega_X \cong \OO_X(\ell)$, and let $G$ be an AG
subscheme of codimension $2$ in $X$.  Then the following
conditions are equivalent:
\begin{itemize}
\item[(i)] There is an  ACM divisor $Y \subseteq X$
satisfying $G_0$ and containing $G$ and an integer $m$ so that $G
\sim mH+M_Y$ on $Y$, where $M_Y$ is the anticanonical divisor.
\item[(ii)] $G$ has an ${\mathcal N}$-type resolution with
${\mathcal N}$ an ACM sheaf of rank $2$ that is an extension of
two rank $1$  ACM sheaves on $X$. In this case we have two exact
sequences:
$$0 \lra \OO_X(\ell-m) \lra \NN \lra \JJ_{G,X} \lra 0,$$
$$0 \lra \JJ_Y \lra \NN \lra \OO_X(Y+\ell-m) \lra 0.$$
\end{itemize}
\end{prop}

\begin{proof} (i) $\Rightarrow$ (ii).  Since $G \sim
mH+M_Y$ on $Y$, we have ${\mathcal J}_{G,Y} \cong \omega_Y(-m)$.
On the other hand, comparing with the ideal sheaf ${\mathcal J}_G$
of $G$ on $X$, we have an exact sequence
\begin{equation}
\label{eq1} 0 \rightarrow {\mathcal J}_Y \rightarrow {\mathcal
J}_G \rightarrow {\mathcal J}_{G,Y} \rightarrow 0.
\end{equation}
We combine this with the natural exact sequence
\begin{equation}
\label{eq2} 0 \rightarrow {\mathcal O}_X \rightarrow {\mathcal
O}(Y) \rightarrow \omega_Y \otimes \omega_X^{\vee} \rightarrow 0
\end{equation}
of \cite[2.10]{Hdivisors}.
%(Here we write ${\mathcal O}(Y)$ for
%the notation ${\mathcal L}(Y)$ of \cite{Hdivisors}.)
Since  $\omega_X \cong {\mathcal O}_X(\ell),$
twisting by $a = \ell - m$ we
get
\begin{equation*}
%\label{eq3}
0 \rightarrow {\mathcal O}_X(a) \rightarrow {\mathcal O}(Y+a)
\rightarrow \omega_Y(-m) \rightarrow 0.
\end{equation*}
Since $\omega_Y(-m) \cong {\mathcal J}_{G,Y}$, we can do the
fibered sum construction with the sequence \eqref{eq1} above and
obtain two short exact sequences
\begin{equation}
\label{eq4} 0 \rightarrow {\mathcal O}_X(a) \rightarrow {\mathcal
N} \rightarrow {\mathcal J}_G \rightarrow 0,
\end{equation}
\begin{equation*}
%\label{eq5}
0 \rightarrow {\mathcal J}_Y \rightarrow {\mathcal N} \rightarrow
{\mathcal O}(Y+a) \rightarrow 0.
\end{equation*}
The first is an ${\mathcal N}$-type resolution of ${\mathcal
J}_G$, and the second shows that ${\mathcal N}$ is an extension of
two rank $1$ ACM sheaves on $X$.

(ii) $\Rightarrow$ (i).  Conversely, suppose given an ${\mathcal
N}$-type resolution of the form \eqref{eq4} above, and suppose
that ${\mathcal N}$ is an extension
\begin{equation*}
%\label{eq6}
0 \rightarrow {\mathcal L} \rightarrow {\mathcal N} \rightarrow
{\mathcal M} \rightarrow 0
\end{equation*}
where ${\mathcal L},{\mathcal M}$ are rank $1$ ACM sheaves on $X$.
The composed map ${\mathcal L} \rightarrow {\mathcal I}_G
\rightarrow {\mathcal O}_X$ shows that ${\mathcal L}$ is
isomorphic to the ideal sheaf ${\mathcal I}_Y$ of an ACM divisor
$Y$ containing $G$.  Then by comparing Chern classes we find
${\mathcal M} \cong {\mathcal O}(Y+a)$.  Dividing the sequence
\eqref{eq4} by ${\mathcal I}_Y$ in the second and third place we
obtain
\begin{equation*}
%\label{eq7}
0 \rightarrow {\mathcal O}_X(a) \rightarrow {\mathcal O}(Y+a)
\rightarrow {\mathcal J}_{G,Y} \rightarrow 0.
\end{equation*}
Comparing with the sequence \eqref{eq2} above, we conclude that
${\mathcal I}_{G,Y} \cong \omega_Y(a-\ell) = \omega_Y(-m)$.
Therefore $G \sim mH+M_Y$ on $Y$.  Note from the isomorphism
${\mathcal I}_{G,Y} \cong \omega_Y(a-\ell)$ that $\omega_Y$ is
locally free at the generic points of $Y$, so $Y$ satisfies $G_0$.
\end{proof}
%\begin{rk} We may ask, when does the extension
%\eqref{eq5} for ${\mathcal N}$ in the above proof split?  It
%splits if and only if there is an effective divisor $Z \sim -Y -
%aH$ such that $C$ is the scheme-theoretic intersection $Y \cap Z$.
%Indeed, if ${\mathcal N}$ splits, then the sequence
%\begin{equation}
%\label{eq8} 0 \rightarrow {\mathcal O}_X(a) \rightarrow {\mathcal
%O}(-Y) \oplus {\mathcal O}(Y+a) \rightarrow {\mathcal I}_C
%\rightarrow 0
%\end{equation}
%gives a map ${\mathcal O}_X \rightarrow {\mathcal O}(-Y-a)$
%defining an effective divisor $Z$, and then the sequence
%\eqref{eq8} shows that $C = Y \cap Z$.  Conversely, if there is
%such a $Z$, then the sequence \eqref{eq8} gives an ${\mathcal
%N}$-type resolution with ${\mathcal N}$ split.
%\end{rk}

From now on $X$ will be a nonsingular cubic surface, $C$ will
denote any nonsingular conic on $X$, $\Gamma$ any twisted cubic on
$X$ and $L$ any of the 27 lines in $X$. The curves $C$, $\Gamma$
and $L$ are arithmetically Cohen-Macaulay curves and, according to
the proof of \cite[2.4]{Hargor}, any arithmetically
Cohen-Macaulay curve on $X$ is linearly equivalent to $C+aH_X$,
$\Gamma+aH$, $L+aH$ or $aH$, for some $a \in \mathbb{N}$. The
sheaves $\JJ_{C,X}$, $\JJ_{\Gamma,X}$ and $\JJ_{L,X}$ are ACM
bundles and their minimal free $\OO_{\PP^3}$-resolutions are:
$$0 \lra \OO_{\PP^3}(-3)^2 \srel{\varphi_1}{\lra} \OO_{\PP^3}(-1)\oplus \OO_{\PP^3}(-2) \lra \JJ_{C,X} \lra 0,$$
$$0 \lra \OO_{\PP^3}(-3)^3 \srel{\varphi_2}{\lra}  \OO_{\PP^3}(-2)^3 \lra \JJ_{\Gamma,X} \lra 0,$$
$$0 \lra \OO_{\PP^3}(-2) \oplus\OO_{\PP^3}(-3) \srel{\varphi_3}{\lra} \OO_{\PP^3}(-1)^2 \lra \JJ_{L,X} \lra 0.$$
Note that $\OO(C)\cong \JJ_L(1)$, $\OO(L) \cong \JJ_{C}(1)$ if $C$
and $L$ are contained in the same plane and $\OO(\Gamma) \cong
\JJ_{\Gamma'}(2)$ if $\Gamma+\Gamma'=2H$.

% Moreover $\JJ_L^{\sigma}
% \cong \JJ_{C}(-1)$ and $\JJ_{C}^{\sigma}\cong \JJ_L(-2)$ if $C +L
% = H$. The first syzygy of  $\JJ_{\Gamma}$ is a rank 2 ACM sheaf
% $\JJ_{\Gamma}^\sigma$. Rank 2 ACM bundles on the cubic surface
% have been classified by Faenzi in \cite{Faenzi}; in his setting
% $\JJ_{\Gamma}$ is a bundle of type (H).
Recall the definition 2.2 of the syzygy  sheaf $\FF^{\sigma}$ of a sheaf $\FF$.

\begin{prop} Let $X \subset \PP^3$ be a nonsingular cubic surface.
\begin{itemize}
\item[a)] If $L$ is a line on $X$, then $\JJ_{L,X}^{\sigma} \cong \JJ_{C}(-1)$ for a conic $C$ in a plane with $L$.
\item[b)] Conversely, $\JJ_{C}^{\sigma}\cong \JJ_L(-2)$.
\item[c)] If $\Gamma$ is a twisted cubic on $X$, then
$\JJ_{\Gamma,X}^{\sigma}$ is a rank 2 ACM sheaf that is an
extension of two ACM line bundles. Indeed, there is an exact
sequence
$$0 \lra \OO_X(L) \lra \JJ_{\Gamma,X}^{\sigma} (3) \lra \OO_X(C) \lra 0$$
where $L$ and $C$ are a line and a conic such that $\Gamma \sim
C+L$.
\end{itemize}
\end{prop}

\begin{proof}
a) and b) are elementary. For c) we proceed as follows. Let $L,C$
be a line and a conic such that $L+C \cong \Gamma$. Then one finds
$C.L=1$, and one can compute $Ext^1(\OO_X(C), \OO_X(L))\cong
H^1(\OO_X(L-C))$ $\cong H^1(\OO_X(L+L'-H))$ where $L'=H-C$. By
duality this is dual to $H^1(\JJ_{L+L',X})$. But $L$ and $L'$ are two skew
lines, so $H^1 \neq 0$ and there exists a non-split extension
$$0 \lra \OO_X(L) \lra \EE \lra \OO_X(C) \lra 0.$$
From this we see that $\EE$ has three global sections, and we can write a diagram
$$\begin{array}{ccccccc}
&  &  & 0 &  & 0 &  \\
&  &  & \downarrow &  & \downarrow &  \\
& 0 & \lra & \FF & \lra & \JJ_{C,X} &  \\
   & \downarrow &    & \downarrow &   & \downarrow &  \\
 0 \lra  & \OO_X & \lra   &  \OO_X^3 & \lra   & \OO_X^2 & \lra 0  \\
& \downarrow &    & \downarrow &   & \downarrow &  \\
0 \lra  & \OO_X(L) & \lra   &  \EE &  \lra   & \OO_X(C) & \lra 0  \\
& \downarrow &    & \downarrow &   & \downarrow &  \\
  & \OO_L(-1) & \lra   &  \GG &  \lra   & 0   \\
&  &    & \downarrow &   &  &  \\
&  &    & 0 &   &  &
\end{array}$$ where $\FF$ and $\GG$ are the
kernel and the cokernel of the map $\OO_X^3 \lra \EE$. By the
snake lemma there is a map $\delta: \JJ_{C,X} \lra \OO_L(-1)$
joining the top and bottom rows into a long exact sequence. If
$\delta=0$ then the extension defining $\EE$  splits, contrary to
hypothesis. Therefore $\delta \neq 0$ and its image must be
$\JJ_{C,X} \otimes \OO_L\cong \OO_L(-1)$ since $L.C=1.$ Thus $\GG=
0$ and we find that $\EE$ is generated by global sections, $\FF
\cong \JJ_{\Gamma,X}$, and so  $\EE^{\sigma} \cong
\JJ_{\Gamma,X}$. This is not quite what we want. But $\EE$ is ACM
and has no free summands because
$Ext^1(\OO_X(a),\JJ_{\Gamma,X})=0$ for all $a \in \ZZ$, so Lemma
\ref{propertiesACM}(vi) applies and we obtain that
$\JJ_{\Gamma,X}^{\sigma} \cong \EE(-3)$.
\end{proof}

We recall the following definition of \cite{CH}.
\begin{defi} An ACM sheaf ${\mathcal E}$ on a normal ACM scheme
$X$ is {\em layered} if there exists a filtration
\[
0 = {\mathcal E}_0 \subseteq {\mathcal E}_1 \subseteq \dots
\subseteq {\mathcal E}_r = {\mathcal E}
\]
whose quotients ${\mathcal E}_i/{\mathcal E}_{i-1}$ are rank $1$
ACM sheaves on $X$ for $i = 1,\dots,r$.
\end{defi}
\begin{prop}\label{syz_layered} Let $X \subset \PP^n$ be a normal arithmetically
Gorenstein scheme. Then the following holds
\begin{itemize}
\item[a)] The dual of a layered sheaf is layered.
\item[b)]  Any extension of layered sheaves is layered.
\item[c)] On a nonsingular cubic surface  $X \subset \PP^3$, a syzygy
of a layered ACM sheaf is stably equivalent to a layered sheaf.
\end{itemize}
\end{prop}

\begin{proof}
a) and b) are obvious. For c) we note that every rank 1 ACM sheaf on $X$ is a twist of $\JJ_{L,X}$, $\JJ_{C,X}$, $\JJ_{\Gamma,X}$ and the syzygy of any of these is layered. The condition then follows from \cite[4.1(d)]{CDH}.
\end{proof}

\begin{ex}\label{ex2}\rm
Take two points $Q,R$ on the cubic surface $X$.  They have an
${\mathcal N}$-type resolution
\[
0 \rightarrow {\mathcal O}_X(-1) \rightarrow {\mathcal N}_{Q+R}
\rightarrow {\mathcal J}_{Q+R} \rightarrow 0.
\]
where $\NN_{Q+R}$ is an ACM bundle by Serre's correspondence. It is
easy to prove that $\NN_{Q+R}$ is indecomposable. Indeed, if
$\NN_{Q+R}$ decomposes as $ \OO(C_1) \oplus \OO(C_2)$, then $C_1$
is an ACM divisor and $C_2 \sim -C_1 -H$. Composing with the map
${\mathcal N}_{Q+R} \rightarrow {\mathcal J}_{Q+R}\hookrightarrow
\OO_X $ we see that $-C_1$ and $C_1+H$ must be effective and both
divisors contain $Q+R$. Numerically, this can only happen if
$-C_1$ is linearly equivalent to a line or a conic in $X$. But as
$Q$ and $R$ are general points in $X$, $Q \cup R$ cannot be
contained in a line or a conic.
%so it is linearly equivalent to $C+aH$ where $C$
%is either a line, a conic or a twisted cubic on $X$, $a$ is an
%integer and $C_2 \sim -C-(a+1)H$.
%As $\NN_{Q+R}(-1)$ has no
%sections we obtain $-2 \leq a \leq 0$ if $C$ is a line or a conic
%and $-3 \leq a \leq 0 $ if $C$ is a twisted cubic. But
%$H^0(\NN_{Q+R})\neq 0$ implies $a=0$ or $-2$ if $C$ is a line or a
%conic and $a=0$ or $a=-3$ if $C$ is a twisted cubic. In any case,
%$c_2(\NN_{Q+R})=-C^2-(2a+1)C.H-a(a+1)H^2$ cannot have degree 2,
%which is a contradiction.

However, the following argument proves that $\NN_{Q+R}$ is
layered. Two general points $Q,R$ are contained in a twisted cubic
curve $\Gamma$, and then $Q+R \sim -K_\Gamma$. Therefore by
\ref{Gorenstein} $\NN_{Q+R}$ is an extension
\begin{equation}\label{extN_2p}
0 \rightarrow {\mathcal O}_X(-\Gamma) \rightarrow \NN_{Q+R}
\rightarrow {\mathcal O}_X(\Gamma-1) \rightarrow 0.
\end{equation}
As $H-\Gamma$ is not effective, this extension does not split.

%On the other hand, if $Q,R$ lie on a conic $C$, we get a similar
%sequence using the conic $C$.  In this case $H-C$ is effective,
%equal to a line $L$. Still, for general $Q,R \in C$, the sequence
%for ${\mathcal N}$ does not split.  Only when $\{Q,R\} = C \cap L$
%does the sequence split.

%Note that this sheaf corresponds to type (E) shifted by -1 in
%Faenzi's classification \cite{Faenzi} and has
Note that $\NN_{Q+R}$ has the following minimal free resolution
over $\OO_X$ of period 2 (see Theorem \ref{eisteo})
\begin{equation}\label{resN2pts}\dots \lra \OO_{X}(-2) \oplus \OO_{X}(-3)^3 \lra
\OO_{X}(-1)^3 \oplus \OO_{X}(-2) \lra \NN_{Q+R} \lra 0.
\end{equation}
\end{ex}
\begin{ex}\label{ex1}\rm
Let $P$ be a point on a nonsingular cubic surface $X$ in ${\mathbb
P}^3$.  Then $P$ has an ${\mathcal N}$-type resolution on $X$
\[
0 \rightarrow {\mathcal O}_X \rightarrow {\mathcal N_P}
\rightarrow {\JJ}_P \rightarrow 0.
\]
One can show that for a general point $P \in X$, the sheaf
${\mathcal N_P}$ is not an extension of ACM line bundles. Indeed,
as $c_1(\NN_P)=0$, such an extension would be of the form $0\lra
\OO(-C) \lra \NN_P \lra \OO(C) \lra 0.$ Composing with the map
${\mathcal N}_{P} \rightarrow {\mathcal J}_{P}\hookrightarrow
\OO_X $ we see that $C$ must be an effective ACM divisor
containing P. As $h^0(\NN_P)=1$, we have $h^0(\OO_X(C))=1$ and $h^0(\OO_X(C-H))=0$, so $C$ must be linearly equivalent
to a line (lines are the only ACM curves $C$ on $X$ for which $h^0(\OO_X(C)=1)$). But a general point $P$ is not contained in a line in
$X$, so such an extension cannot exist. This argument also proves
that $\NN_P$ is indecomposable. Therefore $\NN_P$ is an
indecomposable rank 2 ACM sheaf that is not layered.
%However, if $P$ lies on a line $L$, then $P \sim -H-K_L$ on $L$,
%and ${\mathcal N_P}$ is an extension
%\[
%0 \rightarrow {\mathcal O}(-L) \rightarrow {\mathcal N}
%\rightarrow {\mathcal O}(L) \rightarrow 0.
%\]
%As $-L$ is not effective, the extension does not split.

Note that $\NN_P$
%corresponds to  a sheaf of type (G) in Faenzi's
%classification \cite{Faenzi}. It
has the following 2-periodic minimal free resolution as
$\OO_{X}$-module \begin{equation}\label{resN1pt}\dots \lra
\OO_{X}(-2)^3 \oplus \OO_{X}(-3) \lra \OO_{X} \oplus \OO_{X}(-1)^3
\lra \NN_P \lra 0. \end{equation}

Let us look at the syzygy of $\NN_P$. We consider a complete
intersection of two planes containing $P$  meeting $X$ properly.
Then $P$ is linked to two points $Q \cup R$ by a complete
intersection $Y$ which has the following resolution in $X$
$$0 \lra \OO_X(-2) \lra \OO_X(-1) \oplus \OO_X(-1) \lra \JJ_Y \lra 0 \,.$$
By \cite{CDH} 3.2, $Q\cup R$ has an $\NN$-type resolution as
follows
$$0 \lra \MM \lra \GG \lra \JJ_{Q+R,X} \lra 0$$
where $\GG$ is an extension $0 \lra \OO_X(1)^2 \oplus \OO_X \lra
\GG \lra {\NN_P}^{\sigma\vee} \lra 0$ and $\MM$ is a dissoci\'e
sheaf. As ${\NN_P}^{\sigma\vee}$ is ACM, this extension splits so
$\GG \cong \OO_X(1)^2 \oplus \OO_X \oplus {\NN_P}^{\sigma\vee}.$
On the other hand, $Q \cup R$ has an $\NN$-type resolution as in
example \ref{ex2}, and as any two $\NN$ type resolutions are
stably equivalent, we have $$\OO_X(1)^2 \oplus \OO_X \oplus
{\NN_P}^{\sigma\vee} \oplus \LL_1 \cong \NN_{Q+R}(a) \oplus
\LL_2$$ for some twist $a \in \ZZ$ and some dissoci\'e sheaves
$\LL_1$ and $\LL_2$. But ${\NN_P}^{\sigma\vee}$ and $\NN_{Q+R}(a)$
are indecomposable, so ${\NN_P}^{\sigma\vee}\cong \NN_{Q+R}(a)$.
Looking at the resolutions \eqref{resN2pts} and \eqref{resN1pt} we
find ${\NN_P}^{\sigma\vee}\cong \NN_{Q+R}(2)$. As the first Chern
class of $\NN_{Q+R}$ is $-H$ we obtain ${\NN_P}^{\sigma}\cong
{\NN_{Q+R}}^{\vee}(-2) \cong \NN_{Q+R}(-1).$
%Its syzygy $\NN_P^{\sigma}$ shifted by 2 is a sheaf with second chern class equal to 2 and that corresponds to a bundle of type (E) in Faenzi's classification.

Now by \ref{propertiesACM}(vi), $\NN_P \cong
{\NN_P}^{\sigma\sigma}(3)$ which in turn is isomorphic to
${\NN_{Q+R}}^{\sigma}(2)$. As $\NN_{Q+R}$ is layered, it follows
from \ref{syz_layered} that $\NN_P$ is stably equivalent to a
layered sheaf even though it is not layered itself.
\end{ex}

% Before proving the main theorem of this section we need the
% following easy lemma.
% \begin{lema}\label{lema_filt}
% If $0 \ra \EE' \ra \EE \srel{f}{\ra} \EE'' \ra 0$ is an exact
% sequence of coherent sheaves over a scheme $X$ and $\EE''$ has a
% filtration $0=\FF_0 \subseteq \FF_1 \subseteq  \dots \subseteq
% \FF_{n-1} \subseteq \FF_n=\EE''$ with quotients
% $Q_i=\FF_i/\FF_{i-1}$, then $\EE$ has a filtration $0=\HH_0
% \subseteq \HH_1 \subseteq  \dots \subseteq \HH_{n-1} \subseteq
% \HH_n \subseteq \HH_{n+1}=\EE$ with $\HH_1=\EE'$ and quotients
% $\HH_i/ \HH_{i-1} \cong Q_{i-1}$.
% \end{lema}
% \begin{proof}
% We proceed by induction on $n$. If $n=1$, the result is immediate.
%
% If  $n>1$, We compose $f$ with the map $\pi_n$ is the map $\FF_n=
% \EE'' \lra Q_n$. Then the kernel of $\pi_n f$ is a coherent sheaf
% $\HH_n$ fitting in two exact sequences
% $$0 \lra \HH_n \lra \EE \lra Q_n \lra 0,$$
% $$0 \lra \EE' \lra \HH_n  \lra \FF_{n-1} \lra 0.$$
% We apply induction hypothesis to the last exact sequence and we
% are done.
% \end{proof}

\begin{teo}\label{filtration}
Let $\NN$ be the vector bundle in an $\NN$-type resolution of a
set of general points $Z \subset X$. Then $\NN$ is stably
equivalent to a layered sheaf.
% there is a dissoci\'e sheaf
% $\LL$ and a filtration
% $$0=\FF_0 \subseteq \FF_1 \subseteq  \dots \subseteq \FF_{n-1} \subseteq \FF_n=\NN \oplus \LL$$
% such that the quotients $Q_i=\FF_i/\FF_{i-1}$ are either
% orientable layered rank 2 sheaves or syzygies of orientable
% layered rank 2 sheaves.
\end{teo}
\begin{proof} In \cite{Hargor} the second author proved that any set of general
points $Z$ can be strictly $G$-linked to a general point $P$. We
do induction on the number $t$ of links needed.

If $t=0$, we have one general point $P$ and an $\NN$-type
resolution
$$0 \lra \OO_X \lra \NN_P \lra \JJ_{P,X} \lra 0.$$ Any bundle
$\NN$ corresponding to another $\NN$-type resolution of $P$ is
stably equivalent to $\NN_P$. We have seen in Example \ref{ex1}
that $\NN_P$ is the syzygy of a rank 2 layered sheaf, so this case
is finished.

If $t \geq 1$, there is a strict $G$-link from $Z$ to $Z'$ such
that the induction hypothesis applies to the sheaf $\NN'$
belonging to an $\NN$-type resolution of $Z'$,
$$0 \lra \LL' \lra \NN' \lra \JJ_{Z',X} \lra 0 \,.$$
The strict $G$-link is performed by a Gorenstein scheme $W$ having
an $\NN$-type resolution
$$0 \lra \OO_X(-a) \lra \EE  \lra \JJ_{W,X} \lra 0 $$
where $\EE$ is an extension
\begin{equation}\label{e}0 \lra \OO_X(-Y) \lra \EE \lra \OO_X(Y-aH) \lra 0\end{equation}
for a certain ACM curve $Y \subset X$ and $a \in \ZZ$ (see
\ref{Gorenstein}).
%The curve
%$Y$ is linearly equivalent to either $L+mH$ or $C+mH$ or
%$\Gamma+mH$ for some line $L$, conic $C$ and twisted cubic
%$\Gamma$.
Then by \cite{CDH} Proposition 3.2 we know that there is
an $\NN$-type resolution of $Z$
$$0 \lra \LL \lra \GG \lra \JJ_{Z,X} \lra 0 $$
such that $\GG$ is an extension
$$0 \lra \LL'^{\vee} \oplus \EE^{\vee} \lra \GG \lra \NN'^{\sigma \vee} \lra 0.$$
Any other sheaf $\NN$ appearing in an $\NN$-type resolution $$0
\lra \mathcal{P} \lra \NN \lra \JJ_{Z,X} \lra 0
$$ is stably equivalent to $\GG$, so it is enough to prove that
$\GG$ is stably equivalent to a layered.
 By induction hypothesis there
is a dissoci\'e sheaf $\MM$ such that bundle $\NN' \oplus \MM$ is
layered. By Proposition \ref{syz_layered}c), we have that $(\NN'
\oplus \MM)^{\sigma} \cong \NN'^{\sigma}$ is also stably
equivalent to a layered. Therefore $\NN'^{\sigma \vee}$ is also
stably equivalent to a layered sheaf and so is $\GG$.
\end{proof}

\begin{cor}\label{cor_Ulrichlayered}
Any Ulrich bundle corresponding to a set of general points  is
stably equivalent to a layered sheaf.
\end{cor}

Let $\mathcal{C}$ be the category of ACM bundles on $X$ and  let
$G(\mathcal{C})$ be the corresponding Grothendieck group (i.e  we say that
in $G(\mathcal{C})$ $\EE=\EE`+\EE''$ whenever there is an exact
sequence $0 \lra \EE' \lra \EE \lra \EE''\lra 0$). We regard
$G(\mathcal{C})$ as a $\ZZ[h]$-module with the operation $h \cdot
\EE=\EE(1)$ and define the quotient group
$G'=G(\mathcal{C})/(1-h)G(\mathcal{C})$. In other words, in $G'$
we identify a sheaf with all of its twists.

\begin{cor}
Let $\NN$ be the vector bundle in an $\NN$-type resolution of a
set of general points $Z \subset X$. Then, in $G'$, $\NN$ belongs
to the subgroup generated by rank one ACM bundles and it is
equivalent to $r \OO_X$, $r= \rg(\NN)$.
%there is an $m
%\in \ZZ_+$ such that $\UU_r+m\OO_X$
%belongs to the monoid
%generated by rank one ACM bundles and $\NN_p$ ,
%$\JJ_{\Gamma}^{\sigma}$ where $p$ is a general point on $X$ and
%$\JJ_{\Gamma}^{\sigma}$ are syzygy sheaves of the ideal sheaf of
%one of the 72 twisted cubics on $X$.
%$\{\OO_X, \JJ_L, \JJ_{C_0}, \JJ_{\Gamma}^{\sigma},\NN_p \mid p \textrm{ general in }X$
\end{cor}
\begin{proof}
%We use induction on $r$.
%Set $\UU_1= \NN_P$. In example \ref{ex1and2} we saw that $\NN_P$ is the syzygy of the bundle corresponding to 2 points $\NN_{Q+R}$ and this last bundle is layered. So we have $\NN_P=4\OO_X-\NN_{Q+R}=4\OO_X-\OO(-\Gamma)-\OO(\Gamma)$ for a certain twisted cubic curve as in Example \ref{ex1and2}.

It can be seen in the proof of Theorem \ref{filtration} that, in
$G'$,  $\NN$ is equivalent to
$$\sum_i(\OO(L_i)+\OO(-L_i))+\sum_j(\OO(C_j)+\OO(-C_j))+\sum_k(\OO(\Gamma_k)+\OO(-\Gamma_k))$$
where $L_i$ are lines, $C_j$ are conics and $\Gamma_k$ are twisted
cubics in $X$.

%According to Theorem \ref{filtration}, $\NN$ is stably equivalent
%to a layered sheaf. As $\NN$ is orientable, it remains to prove
%that any orientable layered $\LL$ sheaf is trivial in the group
%$G'$.

%As $\LL$ is layered and orientable, in $G'$ it is equivalent to
%$$\sum_i(\OO(L_i)+\OO(-L_i))+\sum_j(\OO(C_j)+\OO(-C_j))+\sum_k(\OO(\Gamma_k)+\OO(-\Gamma_k))$$
%where $L_i$ are lines, $C_j$ are conics and $\Gamma_k$ are twisted
%cubics in $X$.

%%%%%%%%%%%%%%%%%%%%%%%%%%55
%According to the Theorem and Remark \ref{quotients} $\NN$ admits
%admits a filtration whose quotients are of type $\EE_1$, $\EE_2$,
%$\EE_3$, $\EE_1^{\sigma}$, $\EE_2^{\sigma}$, $\EE_3^{\sigma}$ or
%dissoci\'e sheaves.  We just need to prove that these quotients are
%equivalent to copies of $\OO_X$ in the group $G'$.
%
%$\EE_1$ (resp. $\EE_2$) corresponds to an extension \eqref{e1},
%hence in $G'$ we have $\EE_1=\OO_X(L)+\OO_X(-L)$ (resp.
%$\EE_1=\OO_X(C)+\OO_X(-C)$)
As the syzygy of $\OO_X(-L_i)$ (resp. $\OO_X(-C_j)$) is a twist of
$\OO_X(L_i)$ (resp. $\OO_X(C_j)$), the first two
summands are equivalent to sums of $\OO_X$ in $G'$.
%$\EE_1=2\OO_X$ and
%$\EE_2=2\OO_X$ in $G'$. As a consequence $\EE_1^{\sigma}=2\OO_X$
%and $\EE_2^{\sigma}=2\OO_X.$

%As $\EE_3$ is an extension \eqref{e3}, we have
%$\EE_3=\OO_X(-\Gamma)+\OO_X(\Gamma)$ and
%$\EE_3^{\sigma}=\ell\OO_X-\OO_X(-\Gamma)-\OO_X(\Gamma)$,
%$\ell=2+\rg\EE_3^{\sigma}$.
Let $\Gamma$ be a twisted cubic in $X$. We need to prove that
$\OO_X(-\Gamma)+\OO_X(\Gamma)=2\OO_X$ in $G'$. Take a line $L$
meeting $\Gamma$ in one point. Then there are two conics $D_1$,
$D_2$ such that $L +\Gamma \sim D_1+D_2$ and it is easy to prove
that $(\Gamma-D_1).(\Gamma-D_2)=0.$
% and there is an exact sequence (see
%\cite[6.6]{Faenzi})
%$$0 \lra \OO_X(\Gamma) \lra \OO_X(C_1) \oplus \OO_X(C_2) \lra \OO_X(L) \lra 0.$$
Therefore there is an exact sequence
$$0 \lra \OO_X(\Gamma-L) \lra \OO_X(\Gamma-D_1) \oplus \OO_X(\Gamma-D_2) \lra \OO_X \lra 0$$
and tensoring by $\OO_X(L)$ we obtain the exact sequence
$$0 \lra \OO_X(\Gamma) \lra \OO_X(D_2) \oplus \OO_X(D_1) \lra \OO_X(L) \lra 0.$$
 Therefore, in $G'$ we have
$\OO_X(-\Gamma)+\OO_X(\Gamma)=\OO_X(-D_1)+\OO_X(-D_2)-\OO_X(-L)+\OO_X(D_1)+\OO_X(D_2)-\OO_X(L)$.
Moreover as the syzygy of $\OO_X(-D_i)$  (resp. $\OO_X(L)$) is a
twist of $\OO_X(D_i)$ (resp. of $\OO_X(L)$), we obtain that
$\OO_X(-\Gamma)+\OO_X(\Gamma)=2\OO_X$.
\end{proof}

\bibliographystyle{amsplain}
%amsplain, alpha
%\bibliography{total}

\providecommand{\bysame}{\leavevmode\hbox
to3em{\hrulefill}\thinspace}
\providecommand{\MR}{\relax\ifhmode\unskip\space\fi MR }
% \MRhref is called by the amsart/book/proc definition of \MR.
\providecommand{\MRhref}[2]{%
  \href{http://www.ams.org/mathscinet-getitem?mr=#1}{#2}
} \providecommand{\href}[2]{#2}

\end{document}